\documentclass[12pt]{amsart}
\usepackage{amsmath,amsthm}

\def\pmatrix{\left(\begin{matrix}}
\def\endpmatrix{\end{matrix}\right)}

\def\P{{\mathbb P}}
\def\Z{{\mathbb Z}}
\def\C{{\mathbb C}}

\def\dd{{\rm d}}

\def\de{\delta}

\def\p{\partial}

\def\t{\theta}
\def\s{\sigma}
\def\T{\Theta}
\def\e{\varepsilon}

\def\a{\alpha}
\def\b{\beta}

\def\A{{\mathcal A}}
\def\Airr{\A_g^{irr}}

\def\J{{\mathcal J}}
\def\H{{\mathcal H}}
\def\I{{\mathcal I}}

\def\tt#1#2{{\t\left[\begin{matrix}{#1}\\ {#2}\end{matrix}\right]}}

\theoremstyle{plain}
\newtheorem{thm}{Theorem}

\newtheorem{prop}[thm]{Proposition}

\newtheorem{fact}[thm]{Fact}

\theoremstyle{definition}

\newtheorem{rem}[thm]{Remark}

\def\ee{{\bf e}}
\def\c{{\!\,\cdot\!\,}}
\title{Cubic equations for the hyperelliptic locus}
\author{Samuel Grushevsky}
\address{Mathematics Department, Princeton University, Fine Hall,
Washington Road, Princeton, NJ 08544, USA}
\email{sam@math.princeton.edu}
\thanks{First published in {\it Asian J Math}, {\bf 8} (2004) no.
1, 161-172, special issue dedicated to Yum-Tong Siu on his 60th
birthday. Partially supported by NSF Mathematical Sciences
Postdoctoral Research Fellowship}

\begin{document}
\begin{abstract}
We discuss the conjecture of Buchstaber and Krichever from
\cite{krbu} that the multi-dimensional vector addition formula for
Baker-Akhiezer functions obtained there characterizes Jacobians
among principally polarized abelian varieties, and prove that it
is indeed a weak characterization, i.e. that it is true up to
additional components. We also show that this addition formula is
equivalent to Gunning's multisecant formula for the Kummer variety
obtained in \cite{gu2}.

We then use the computation of the coefficients in the addition
formula from \cite{krbu} to obtain cubic relations among theta
functions that (weakly) characterize the locus of hyperelliptic
Jacobians among irreducible abelian varieties. In genus 3 our
equations are equivalent to the vanishing of one theta-null, and
thus are classical (see \cite{mumford}, \cite{poor}), but already
for genus 4 they appear to be new.
\end{abstract}
\maketitle

\section{Definitions and notations}
We work over $\C$, and fix the dimension/genus $g>1$. Let $\H_g$
be the Siegel upper half-space --- the set of all $g\times g$
period matrices $\tau$, i.e. symmetric complex $g\times g$
matrices with positive definite imaginary part. Each such $\tau$
corresponds to an abelian variety $X_\tau:=\C^g/\tau \Z^g+\Z^g$,
and the moduli space $\A_g$ of principally polarized abelian
varieties (ppavs) is then the quotient of $\H_g$ by a certain
action of the symplectic group $\operatorname{Sp}(2g,\Z)$.

A ppav is called irreducible if it is not isomorphic to a product
of two lower-dimensional ppavs (with polarization). For
convenience we denote by $\Airr$ the moduli space of irreducible
ppavs of genus $g$. When in the following we say ``abelian
variety'', we actually mean a ppav.

Denoting $\ee(x):=\exp(\pi i x)$, for a period matrix $\tau$ and a
vector $z\in\C^g$ we define the theta function with
characteristics $\e,\de\in(\Z/2\Z)^g$, thought of as vectors
consisting of zeros and ones, to be
$$
\tt\e\de(\tau,z):=\sum\limits_{m\in\Z^g} \ee \left[\left(
m+\frac{\e}{2},\tau(m+\frac{\e}{2})\right)+2\left(m+\frac{\e}{2},z+
\frac{\de}{2}\right)\right],
$$
where $(\cdot,\cdot)$ denotes the scalar product. A theta function
with characteristics is even or odd as a function of $z$ depending
on whether the scalar product $(\e,\de)$ is even or odd.

We denote by $\t(\tau,z):=\tt00(\tau,z)$ the classical Riemann's
theta function. Theta functions with characteristics are, up to a
constant factor, just the values of Riemann's theta function of a
shifted argument:
\begin{equation}
\label{shift}
\begin{array}{l}
\t(\tau,z+\frac{\tau\e+\de}{2})=
\sum\limits_{m\in\Z^g}\ee \left[(m,\tau m)+2\left(m,z+\frac{\tau
\e+\de}{2}\right)\right]\\
\quad =\!\sum\limits_{m\in\Z^g}\!\ee\left[\left(
m\!+\!\frac{\e}{2},\tau(m\!+\!\frac{\e}{2})\right)
+2\left(m\!+\!\frac{\e}{2},z\!+\!\frac{\de}{2}\right)
-(\frac{\e}{2},\tau\frac{\e}{2})-\left(\e,z+\frac{\de}{2}\right)\right]\\
\quad
=(-1)^{(\e,\de)}\ee\left[-\frac{1}{4}(\e,\tau\e)-(\e,z)\right]\tt\e\de(\tau,z).
\end{array}
\end{equation}
Thus instead of thinking of a characteristic $\left[\begin{matrix}
\e\\ \de\end{matrix}\right]$ as two integer vectors it sometimes
is better to think of it as the point $\frac{\tau\e+\de}{2}$ of
order two on the abelian variety $X_\tau$.

We further define theta functions of the second order to be
$$
\T[\e](\tau,z):=\tt{\e}{0}(2\tau,2z).
$$
For a fixed $\tau$ the theta functions, as
functions of the variable $z$, are sections of certain bundles on
the abelian variety $X_\tau$, which is to say that if the variable
$z$ is translated by a vector of the lattice $\tau\Z^g+\Z^g$,
theta functions multiply by a certain number. In fact it is known
that all theta functions of the second order are sections of the
same bundle, denoted $2\T$, and transform as follows:
\begin{equation}
\label{autom} \T[\e](\tau, z+e_j+\tau e_k)=\ee(-2(e_k,\tau
e_k)-4(e_k,z))\T[\e](\tau,z),
\end{equation}
where we denote by $e_k$ the basis vector for the $k$'th direction
in $\C^g$.

Theta functions of the second order form a basis for the sections
of $2\T$ over $X_\tau$. The square of any theta function with
characteristics is also a section of $2\T$, and thus is
expressible as a linear combination of theta functions of the
second order. In fact a slightly more general formula, Riemann's
bilinear addition theorem, holds:
\begin{equation}
\label{tT}
\tt\e\de(z)\tt\e\de(w)=\sum\limits_{\s\in(\Z/2\Z)^g}(-1)^{(\de,\s)}
\T[\s+\e]\left(\frac{z+w}{2}\right)\T[\s]\left(\frac{z-w}{2}\right).
\end{equation}

For a fixed $\tau$ the map $z\to\lbrace\T[\e](\tau,z)\rbrace_{\rm
all\ \e}$ defines the Kummer embedding $K:X_\tau/\pm
1\to\P^{2^g-1}$. This map is well-defined since all theta
functions of the second order are sections of the same line
bundle, and are even in $z$.

\smallskip
The values of theta functions at $z=0$ are called the associated
theta constants. Theta constants of the second order are modular
forms of weight one half with respect to a certain finite index
normal subgroup $\Gamma(2,4)\subset \operatorname{Sp}(2g,\Z)$,
which is to say that if we act upon $\tau$ by some
$\gamma\in\Gamma(2,4)$, then $\T[\e](\gamma\tau)=k(\gamma,\tau)
\T[\e] (\tau)$, where $k$ is some multiplier depending on $\gamma$
and $\tau$, but independent of $\e$. Thus letting
$\A_g(2,4):=\H_g/ \Gamma(2,4)$, we see that theta constants of the
second order define a map $Th:\A_g(2,4)\to \P^{2^g-1}$, which is
known to be generically injective for all genera, and injective
for $g\le 3$. The level moduli space $\A_g(2,4)$ is a finite cover
of $\A_g$.

Let us denote by $\J_g\subset\A_g$ the locus of Jacobians of
Riemann surfaces of genus $g$, and by $\I_g\subset\J_g$ the locus
of Jacobians of hyperelliptic Riemann surfaces. The question of
characterizing $\J_g$ within $\A_g$ is called the Schottky
problem, and that of characterizing $\I_g$ --- the Schottky
problem for the hyperelliptics. More precisely, one takes the
preimages $\J_g(2,4)$ and $\I_g(2,4)$ of $\J_g$ and $\I_g$,
respectively, under the covering map $\pi:\A_g(2,4)\to\A_g$, and
asks to describe $Th(\J_g(2,4))$ and $Th(\I_g(2,4))$ inside
$Th(\A_g(2,4))$. The question of describing
$Th(\A_g(2,4))\subset\P^{2^g-1}$, i.e. determining all the
relations in the subring of the ring of modular forms generated by
theta constants is also of interest, but we will not discuss it
here. Notice that $\J_g$ and $\I_g$ are irreducible, while
$\J_g(2,4)$ and $\I_g(2,4)$ have many irreducible components. We
refer the reader to \cite{ig} for more details on theta functions,
and to \cite{gr2} for more details on the Schottky problem.

We will always think of a curve $C$ embedded in its Jacobian by
the Abel-Jacobi map $A:C\to J(C)$ with some choice of the basis
for the space of holomorphic differentials and of the starting
point $P$ made. This choice will be made explicitly when
necessary. To avoid technical difficulties in the following
sections, it will often be easier to work with the universal cover
$\tilde C$ of a curve $C$ and the universal cover $\C^g$ of the
abelian variety $X_\tau$, and later take the automorphy properties
of theta functions into account. The abelian variety will be fixed
throughout, and thus we will often omit $\tau$ from the notations
for theta functions and constants.

\smallskip In this work we first prove the conjecture of Buchstaber and
Krichever stated in \cite{krbu} that the validity of a certain
$g$-dimensional addition formula developed there and in
\cite{krbu1} characterizes Jacobians, {\it but only up to
additional components} (i.e. that the locus of Jacobians is an
irreducible component of the locus where the addition formula is
satisfied), and then proceed to obtain from this some explicit
identities for theta functions of hyperelliptic curves, using the
explicit coefficients for the addition formula from \cite{krbu}.

In \cite{mumford}, theorem 9.1, and references therein Mumford
showed that the hyperelliptic locus is characterized by a certain
set of vanishing and non-vanishing conditions for theta constants
with characteristics (the idea goes back at least to Thomae,see
 \cite{thomae}). In \cite{poor} Poor showed that on $\Airr$
Mumford's vanishing conditions by themselves (without the
non-vanishing) define precisely the hyperelliptic locus, i.e. that
there are no extra components. However, it is still not known how
to obtain an ideal-theoretic description of the closure of
$Th(\I_g)$ inside $Th(\A_g)$. It is known that if the vanishing
holds and we also have some vanishing instead of non-vanishing,
the abelian variety must be reducible, but then it does not
necessarily have to be a limit of hyperelliptic Jacobians. Thus it
would be interesting to study our equations on the reducible
locus. It would also be very interesting to compare Mumford's
equations to ours, but we have not been able to achieve this yet.

We would also like to refer to a recent work \cite{sm} for a
further discussion of these issues as well as a description of
components of $\I_g^{2,4}$ as locally complete intersections.

\section*{Acknowledgements}
We would like to thank Igor Krichever for bringing to our
attention and explaining to us \cite{krbu1} and \cite{krbu} and
the conjectural characterization of Jacobians by the addition
formula, which got this work started, Robert Gunning for helpful
discussions on addition formulas for theta functions and
multisecants, and Emma Previato for valuable discussions of the
equations for Kummer and modular varieties.

We would like to especially thank Giuseppe Pareschi and Mihnea
Popa, the authors of related recent preprint \cite{popa} (see
remark \ref{pprem} for details), for pointing out that the
assumption of theorem 3 in the published version of this paper was
too weak, and thus that some general position assumption was
needed. The current text is the corrected version of the paper
(the mistake was that lemma 4 was false), and the appropriate
erratum is to appear in print shortly.

\smallskip
I owe a debt of gratitude to Professor Yum-Tong Siu for teaching
me and sharing with me his insights and ideas in moduli theory. It
is an honor for me to dedicate this paper to Professor Yum-Tong
Siu on the occasion of his 60th birthday.

\section{Addition formula and multisecants}
From \cite{krbu} we know the following addition formula (called
``formula'' to distinguish it from Riemann's addition
``theorem''):
\begin{fact}[\cite{krbu}, theorem 1]
Let $P:=A_0,A_1,\ldots,A_g,Q:=A_{g+1}\in\tilde C\subset\C^g$, and
denote by $R:=K-\sum\limits_{i=1}^g A_i$ the vector of Riemann
constants shifted by $-\sum A_i$. Then for all $x,y\in\C^g$ the
following identity is satisfied:
$$
\frac{\t(Q+x+y+R)\t(P+R)}{\t(Q+R)\t(P+x+y+R)}=
$$
$$
=\frac{\t(Q+x+R)\t(P+R)}{\t(Q+R)\t(P+x+R)}\cdot\frac{\t(Q+y+R)\t(P+R)}
{\t(Q+R)\t(P+y+R)}
$$
$$
-\sum\limits_{k=1}^g\frac{\t(Q+R+A_k-P)\t(Q+R-A_k+P+x+y)}{\t^2(Q+R)
\t(P+x+y+R)\t(2A_k-P+R)}
$$
$$
\cdot\frac{\t(R+A_k+x)\t(P+R)}{\t(P+x+R)}\cdot\frac{\t(R+A_k+y)\t(P+R)}{\t(P+y+R)}.
$$
\end{fact}

Though this formula may look formidable, it is very explicit and
is written entirely in terms of theta functions. In the following,
we take the Abel-Jacobi map to start at $P$, so that $P=0\in\C^g$.
Upon cancellations and multiplication by the common denominators,
the above formula becomes simply
$$
\t(Q+x+y+R)\t(Q+R)\t(x+R)\t(y+R)
$$
$$
=\t(x+y+R)\t(R)\t(Q+x+R)\t(Q+y+R)
$$
$$
-\sum\limits_{k=1}^g\frac{\t(R)\t(Q\!+\!A_k\!+\!R)\t(Q\!-\!A_k\!+\!
x\!+\!y\!+\!R)\t(A_k\!+\!x\!+\!R)\t(A_k\!+\!y\!+\!R)}{\t(2A_k\!+\!R)}.
$$
To see that this is in fact equivalent to Gunning's general
multisecant formula from \cite{gu2} (see Poor's work \cite{poor}
for an in-depth discussion) we use Riemann's bilinear addition
theorem for the last two factors of each term. Denoting
$z:=\frac{x+y}{2}$ and $w:=\frac{x-y}{2}$, notice that the
half-difference is always simply $\frac{x-y}{2}$, so we will have
a common factor of $\T[\s](w)$, and the resulting equation will be
$$
\quad
0=\!\!\sum\limits_\s\!\Big[\t(Q\!+\!2z\!+\!R)\t(Q\!+\!R)\T[\s](z\!+\!R)
-\t(2z\!+\!R)\t(R)\T[\s](Q\!+\!z\!+\!R)
$$
\begin{equation}
\label{mess}
+\!\!\sum\limits_{k=1}^g\!\frac{\t(R)\t(Q\!+\!\!A_k\!+\!R)\t(Q\!-\!
\!A_k\!+\!2z\!+\!R)\T[\s](\!A_k\!+\!z\!+\!R)}{\t(2A_k\!+\!R)}\Big]
\T[\s](w)
\end{equation}

In the above the coefficient in the square brackets does not
depend on $w$, so we have an equation $\sum_\s b_\s\T[\s](w)=0\
\forall w$, where the coefficients $b_\s$ do not depend on $w$.
Since theta functions of the second order are a basis for sections
of $2\T$ and thus linearly independent, it means that all
coefficients $b_\s$ must be zero. Then since $b_\s$ actually are
some functions of $z$ multiplied by $\T[\s](A_i+z+R)$ ($i$ here
ranges from $0$ to $g+1$, i.e. includes $P$ and $Q$), we have
$$
\forall \s,\forall z\quad \sum\limits_{i=0}^{g+1}
c_i(z,A_i)\T[\s](A_i+z+R)=0
$$
for appropriate $c_i$'s.

Since $R$ does not depend on $z$, in the above we can shift $z$ by
$R$ and redefine the $c_i$ to see that the addition formula of
\cite{krbu} implies the existence for any $z$ of some complex
numbers $c_i$, not all simultaneously zero, such that for some
fixed $A$'s lying on the image $A(C)\subset J(C)$ we have
\begin{equation}
\label{collinear} \sum\limits_{i=0}^{g+1}c_i(z)K(A_i+z)=0
\qquad\forall z\in\C^g,
\end{equation}
which is equivalent to saying that the $g+2$ points $K(A_i+z)$ lie
on a $g$-plane in $\P^{2^g-1}$. In fact Gunning proves a more
general theorem (with a different and seemingly much more
complicated expression for $c_i$'s in terms of the prime form):

\begin{fact}[\cite{gu2}, theorem 2] For any curve $C$ of genus $g$,
for any $1\le m\le g$ and for any points $x_1,\ldots,x_m,A_0,
\ldots, A_{m+1}\in A(C)$ the $m+2$ points $K(A_i+\sum x-\sum A)$
are collinear, i.e. lie on the intersection of the Kummer variety
of $C$ with an $m$-plane in the projective space $\P^{2^g-1}$.
\end{fact}
In particular since the $g$'th symmetric power of the curve is its
Jacobian, $S^g(A(C))=J(C)$, the case $m=g$ of this theorem is
formula (\ref{collinear}), while the case $m=1$ is the case of a
family of trisecant lines. It is shown in \cite{gu1} that the
existence of a family of trisecants characterizes Jacobians among
irreducible ppavs.

\section{Characterizing Jacobians by families of multisecants}
We will now show that Buchstaber-Krichever's addition formula and
Gunning's multisecant formula weakly characterize Jacobians.
\begin{thm}
Let $X$ be an irreducible principally polarized abelian variety of
dimension $g$, and let $A_0, \ldots, A_{g+1}$ be distinct points
of $X$. Suppose that $\forall z\in X$ the $g+2$ points $K(A_i+z)$
are linearly dependent. Assume moreover the following general
position condition: that there exist some $k$ and $l$ such that
for $y:=-\frac{A_k+A_l}2$ the linear span of the points $K(A_i+y)$
for $i=0\ldots g+1$ is of dimension precisely $g+1$, and not less.
Then $X$ is the Jacobian of some curve $C$, and all $A_i\in A(C)$.
\end{thm}
\begin{proof}
Assume $g\ge 4$, otherwise the theorem is trivial. Working in the
spirit of Gunning's work \cite{gu1}, we reduce the theorem to the
case of the trisecant.

Since the rank of the $(g+2)\times 2^g$ matrix $K(A_i+y)$ is equal
to exactly $g+1$, it means that there exists unique $\vec
c(y)\in\P^{g+1}$ such that $\sum c_i(y)K(A_i+y)=0$. Moreover,
since we have $A_k+y=-(A_l+y)$ and thus $K(A_k+y)=K(A_l+y)$, we
must then have $p:=\vec c(y)={\bf 1}_k-{\bf 1}_l$, where ${\bf
1}_i$ denotes the basis vector of $\C^g$ in the $i$'th direction.
Since the rank of $K(A_i+y)$ is equal to $g+1$, the rank of
$K(A_i+z)$ must be equal to $g+1$ identically in a neighborhood of
$y$, and thus for all $z$ sufficiently close to $y$ there is a
unique projective solution $\vec c(z)$ to $\sum c_i(z)K(A_i+z)=0$.

We will show that the differential $\dd \vec c:T_y\C^g\to
T_p\P^{g+1}$ is of maximal rank. This will then imply that locally
near $p$ the image $\vec c(\C^g)$ is of dimension $g$ and thus
locally the preimage of any coordinate plane $\P^2\subset\P^{g+1}$
containing $p$ is at least one-dimensional. Then we have a
one-dimensional family of trisecants of the Kummer variety, and by
Welters' \cite{welters} infinitesimal version of Gunning's
trisecant criterion from \cite{gu1}, $X$ is a Jacobian of some
curve $C$ with the points $A_i,A_k,A_l$ lying on $A(C)$. Thus we
see that all $A_i$ lie on $A(C)$.

\smallskip
For contradiction, suppose that the rank of $\dd \vec c|_y$ is not
maximal, i.e. that there is some vector $v\in \C^g$ such that we
have $\frac{\p}{\p v} c_i(z)|_y=\lambda c_i(y)$ for some constant
$\lambda$ independent of $i$ --- this means that the derivative
$\frac{\p}{\p v}$ of the projective point $\vec c$ is zero.

Let us take the derivative $\frac{\p}{\p v}$ of (\ref{collinear})
at $y$:
$$
\begin{matrix}
\sum\limits_{i=0}^{g+1}\left.\frac{\p c_i(z)}{\p v}
K(A_i+z)+c_i(z)\frac{\p K}{\p v}(A_i+z)\right|_y\\
=\lambda(K(A_k+y)-K(A_l+y))+(
\frac{\p K}{\p v}(A_k+y)-\frac{\p K}{\p v}(A_l+y))\\
=2 \frac{\p K}{\p v}(A_k+y)=0,
\end{matrix}
$$
because $A_k+y=-(A_l+y)$, theta functions are even and their
derivatives are odd. But this then implies that $\frac{\p K}{\p
v}(A_k+y)=0$, which is impossible by Wirtinger's theorem unless
$A_k+y=\frac{A_k-A_l}{2}$ is a point of order two in $X$. If this
is the case, though, it would mean that $A_k-A_l=0\in X$, which
contradicts the assumption that all $A_i$ are distinct. Thus we
have arrived at a contradiction and showed that the differential
$\dd\vec c: T_y\C^g\to T_P\P^{g+1}$ is injective.
\end{proof}

\begin{rem}\label{pprem} In their recent preprint \cite{popa}
Pareschi and Popa prove the following statement.

\smallskip
\noindent {\bf Castelnuovo-Schottky lemma (\cite{popa}).} Let
$A_0\ldots A_n\in X$ be a set of points in $\t$-general position
imposing only $g+1$ conditions on $|2\T|_a$ for $a$ general. Then
if $n\ge g+1$, $X$ is the Jacobian of some curve $C$, and moreover
$A_i\in A(C)$.

\smallskip
All sections of $|2\T|_a$ are scalar products $v\cdot K(z+a)$ for
some $2^g$-dimensional vector $v$, and the condition for a section
to vanish at $A_i$ is simply that $v\cdot K(A_i+a)=0$. Thus the
points imposing only $g+1$ conditions on $|2\T|_a$ means that for
general $a$ the rank of the $(n+1)\times 2^g$ matrix $K(A_i+a)$ is
equal to $g+1$, which is the collinearity condition of our theorem
3.

The $\t$-general position means that for any $g+1$ points
$p_1\ldots p_g,p_{g+1}$ among $\lbrace A_i\rbrace$ there exists
some translate $a\in X$ such that $\t(p_i+a)=0\ne \t(q+a)$ ---
this is different from our rank condition. Since both our theorem
3 and the Castelnuovo-Schottky lemma characterize Jacobians, these
assumptions must be equivalent, but we do not know a way to derive
one from the other directly.
\end{rem}

With more work it can also be shown that Gunning's addition
formula for all other values of $m$ also serves to characterize
Jacobians:

\begin{prop}
Let $X$ be an irreducible principally polarized abelian variety of
dimension $g$, and let $A_0, \ldots, A_{m+1}$ be different points
of $X$. Suppose that the $m+2$ points $K(A_i+z)$ are linearly
dependent for any $z\in M$, where the set $M\subset X$ is at least
$m$-dimensional at the point $y=-\frac{A_k+A_k}{2}$ for some $k$
and $jk$ (i.e. has a non-degenerate $m$-jet at $y$). Assume
moreover that the rank of the $(m+2)\times 2^g$ matrix $K(A_i+y)$
is equal to precisely $m+1$. Then $X$ is the Jacobian of some
curve.
\end{prop}
\begin{proof}
We again study the local situation near $y$ and work with the germ
of $M$ at this point. Then we want to show that $\dd\vec c$ is
injective on $T_y\tilde M$. But since $\dd \vec c$ is
non-degenerate on $T_y\C^g$, it is also non-degenerate on a
subspace, and we are done.
\end{proof}
\begin{rem}
By imitating Gunning's proof of the trisecant theorem in
\cite{gu1} it seems to be also possible to show that if $M$ is
$m$-dimensional at some point, not necessarily
$-\frac{A_i+A_j}{2}$ and under some general position assumption,
then $X$ a Jacobian. However, as we do not need this result in
what follows, and the proof would be technically rather
complicated, we will not give it.
\end{rem}

In \cite{gu2} Gunning shows that for a Jacobian $X$ with fixed
$A$'s the coefficients $c_i$ in (\ref{collinear}) are unique up to
scaling and expressible in terms of the Klein-Gunning prime form.
The Klein-Gunning prime form expression is rather hard to deal
with (see \cite{poor} for a detailed discussion and computations);
however, the expression for $c_i$ obtained in formula (\ref{mess})
seems more amenable.

Equation (\ref{mess}) includes, however, both theta functions of
the second order and the classical Riemann's theta function with
zero characteristic. Let us use the addition theorem once again
for the last two factors of the type $\t(\ldots)\t(\ldots)$ in
each term. We then get for all $\s$
\begin{equation}
\label{lastadd}
\begin{matrix}
\sum\limits_\e\!\T[\e](Q\!+\!z\!+\!R)\T[\e](z)\T[\s](z\!+\!R)
-\T[\e](z\!+\!R)\T[\e](z)\T[\s](Q\!+\!z\!+\!R)\\
+\sum\limits_{k,\e}\frac{\t(R)}{\t(2A_k+R)}\T[\e](Q+z+R)\T[\e]
(z-A_k)\T[\s](A_k+z+R)=0.
\end{matrix}
\end{equation}
Since all the transformations that we have done so far are
equivalencies, we do not lose any information, unless {\it all}
the coefficients $c_i$ in the formula above are identically zero
--- in which case the projective point $\vec c(z)$ is never defined.

Since the coefficients in (\ref{collinear}) are unique, they must
be exactly the ones given by formula (\ref{lastadd}), and thus
formula (\ref{lastadd}) is satisfied if and only if the
collinearity condition (\ref{collinear}) is satisfied.

\begin{prop} With the same assumptions and notations as above,
identity (\ref{lastadd}) characterizes Jacobians among all
irreducible abelian varieties.
\end{prop}

\section{Addition formula for the hyperelliptic case}
Formula (\ref{lastadd}) we obtained is in terms of theta functions
evaluated at different points. In the classical approach to the
Schottky problem (see \cite{sch} for the origins and
\cite{farkas}, \cite{gr2} for a review) one wants to characterize
the Jacobian locus by some algebraic relations among theta
constants. Thus it would be nice if for some special values of
$z,P,Q$ and $A$'s the addition formula yielded such equations.

However, from the transformation rule (\ref{autom}) for theta
functions it is easy to see that if a vector $v$ is not a point of
order two, then $\T[\e](z+v)$ is not a section of the bundle $2\T$
and thus cannot be expressed as a linear combination of theta
functions of the second order. Thus we only have a reasonable hope
of getting from (\ref{lastadd}) some equations for theta constants
if we are so lucky that all the ``shifts'' of $z$ that appear
there are points of order two. In particular, this means that the
points $A_i+z+R$ must be of order two for all $i$.

Now suppose that indeed both $A_i+z+R$ and $A_j+z+R$ for $i\ne j$
are points of order two on $J(C)$. Then their difference,
$A_i-A_j$, is also of order two, so we have $2A_i-2A_j=0\in J(C)$.
By Abel's theorem this then means that there is a function $f$ on
$C$ whose divisor is equal to $2A_i-2A_j$, i.e. with a double pole
at $A_j$ and holomorphic on $C-\lbrace A_j\rbrace$. Since the
existence of such a function characterizes hyperelliptic curves,
it means that $C$ then has to be hyperelliptic. For the
hyperelliptic curves it is known (see, for example,
\cite{mumford}) that if we take $P=0$ to be the image of one of
the Weierstrass points, then the other $2g+1$ Weierstrass points
will also map to points of order two on the Jacobian.

Thus let us assume that all $A$'s and $Q$ in formula \ref{lastadd}
are chosen to be points of order two, i.e. that we are dealing
with a hyperelliptic curve, and let us rewrite the addition
formula in this case. Denote $Q=\frac{\tau\a_0+\b_0}{2}$,
$A_k=\frac{\tau\a_k+\b_k}{2}$ and $R=\frac{\tau\a+\b}{2}$ --- in
fact $R$ is expressible in terms of $A$'s and Riemann constants,
but we will deal with this later.

Now we rewrite formula (\ref{lastadd}) for these $A$'s and $Q$. In
doing this, we need to be extra careful to remember that we are
actually working on $\tilde C$ and $\C^g$, as not to omit any
important automorphy factors. Indeed, from the automorphy
properties of $\t$ it follows that
$$
\frac{\t(R)}{\t(2A_k+R)}=\frac{\t\left(\frac{\tau\a+\b}{2}\right)}
{\t\left(\tau\a_k+\b_k+\frac{\tau\a+\b}{2}\right)}=\ee[(\a_k,\tau
\a_k)+(\a,\tau\a_k)].
$$
Also for any integers $a$ and $b$ it follows from (\ref{shift})
and (\ref{autom}) that
\begin{equation}
\begin{array}{l}
\label{change} \T[\de]\left(\tau,z+\frac{\tau a+b}{2}\right)=
\tt\de0(2\tau,2z+\tau a+b)\\
\quad\quad =\t(2\tau,2z+\tau\de+\tau a)\ee\left[\frac12(\de,\tau\de)
+2(\de,z)+(\de,\tau a)+(\de,b)\right]\\
\quad\quad=\tt{\de+a}0(2\tau,2z)\ee\left[-\frac12(\de+a,\tau(\de+a))
-2(\de+a,z)\right]\\
\quad\quad\quad\cdot\ee\left[\frac12(\de,\tau\de) +2(\de,z)
+(\de,\tau a)+(\de,b)\right]\\
\quad\quad =(-1)^{(\de,b)}\ee\left[-\frac12(a,\tau
a)-2(a,z)\right]\T[\de+a](\tau,z).
\end{array}
\end{equation}
When we substitute this into (\ref{lastadd}) notice that as
functions of $z$ all terms are actually sections of the same
bundle, $6\Theta$, as each is cubic in theta functions of the
second order. Thus the $\ee[(*,z)]$ factors must cancel
everywhere. Also evaluating at $z=0$ and noticing that all terms
are modular forms in $\tau$ with respect to $\Gamma(2,4)$ of the
same weight, we expect the factors $\ee[(*,\tau*)]$ to cancel as
well. A trivial but tedious computation confirms this, and we
arrive at

\begin{prop}
An irreducible abelian variety $X_\tau$ with some points
$P=0,Q,A_1,\ldots, A_g$ with $A_i=\frac{\tau \a_k+\b_k}{2}\in
\C^g$ is the Jacobian of a hyperelliptic curve $C$, and
$P,Q,A_i\in A(\tilde C)\subset \C^g$ if and only if the following
is satisfied for all $\s\in(\Z/2\Z)^g$ and for all $z\in \C^g$:
\begin{equation}
\label{addhyp}
\begin{array}{cl}
&\sum\limits_\e(-1)^{(\e,\b+\b_0)}\T[\e+\a+\a_0](z)
\T[\e](z)\T[\s+\a](z)\\
=&\sum\limits_{e,k}(-1)^{(\e,\b+\b_0+\b_k)+(\s,\b_k)}
\T[\e\!+\!\a\!+\!\a_0](z)\T[\e\!+\!\a_k](z)\T[\s\!+\!\a\!+\!\a_k](z)\\
+&\sum\limits_\e(-1)^{(\e,\b)+(\s,\b_0)}\T[\e+\a](z)\T[\e](z)\T[\s+\a+\a_0](z).
\end{array}
\end{equation}
and not all the coefficients in front of the Kummer images (i.e.
in front of the last theta function factor, for $\s$ varying)
appearing here are identically zero in $z$.
\end{prop}

\section{Cubic equations for the hyperelliptic locus}
To make formula (\ref{addhyp}) entirely explicit, we now need to
pick some specific way to map a hyperelliptic curve into its
Jacobian, and pick some $g+2$ Weierstrass points on it in a
certain way. This is indeed a very classical construction.

Let us think of a hyperelliptic curve sitting on a skewer that
intersects it in precisely the $2g+2$ Weierstrass points. Label
them $p_1,\ldots, p_{2g+2}$ going from left to right along the
skewer. Then pick for the basis of the cycles $a_i$ to be the loop
around the $i$'th handle, passing through points $p_{2i-1}$ and
$p_{2i}$, and $b_i$ to be the loop around the $i$'th hole, passing
through $p_{2i}$ and $p_{2i+1}$. Then thinking of the skewer as
being the $x$ axis and the whole picture being that of $y^2=\prod
(x-p_i)$, we can compute the images of $p_i$ in the Jacobian.
Indeed, let us use $p_1$ as the starting point, so that
$A(p_1)=0$. Then we see that $A(p_2)=\frac{e_1}{2}$, $A(p_3)=
\frac{\tau e_1+e_1}{2}$, $A(p_4)=\frac{\tau
e_1+e_1}{2}+\frac{e_1+e_2}{2}=\frac{\tau e_1+e_2}{2}+e_1$, $A(p_5)
=\frac{\tau (e_1+e_2)+e_2}{2}+e_1$, and in general we have
$A(p_{2i})=\frac{\tau s_{i-1}+e_i}{2}+s_{i-1}$ and
$A(p_{2i+1})=\frac{\tau s_i+e_i}{2}+s_{i-1}$ for $1\le i\le g$,
while $A(p_{2g+2})=\frac{\tau s_g}{2}$ (where for convenience we
have denoted $s_k:=\sum\limits_{i=1}^k e_i$).

For our purposes let us choose $P:=0=A(p_1)$, $Q:=A(p_2)$ and
$A_k:=A(p_{2k+2})$ for $1\le k\le g$. We now need to compute the
vector $R$, i.e to compute the vector of Riemann constants and
subtract from it the sum of $A$'s. The result is certainly
classical: $R=Q=A(p_2)$. To prove this one can note that by
definition $R$ is the unique vector such that $\t(A(p)+R)$ as a
function of $p\in C$ has precisely $g$ zeroes at $A^{-1}(A_k)$,
i.e. at $p_{2k+2}$. To check that this is the case we note that
$A(p_2)+A(p_{2i})$ is always odd, as a theta characteristic, so
that even Riemann's theta function will vanish at the point
$A(p_2)+A(p_{2i})$, while $A(p_2)+A(P_{2i+1})$ is even; thus
$R=A(p_2)$.

Let us now substitute all this into formula (\ref{addhyp}). We
have $\a=\a_0=0$, $\b=\b_0=e_1$, $\a_k=s_k$ and $\b_k=e_{k+1}$,
where we understand $e_{g+1}$ to be zero. Since all theta
functions of the second order are periodic with respect to $z\to
z+e_i$, the additional integer shifts by $s_i$ do not matter, and
finally (\ref{addhyp}) yields
\begin{thm}
\label{thmchar} An irreducible period matrix $\tau\in\H_g$ is the
period matrix of a hyperelliptic Jacobian with the basis of cycles
chosen as above if and only if the following cubic identity for
theta functions of the second order is satisfied for all
$\s\in(\Z/2\Z)^g$ and for all $z\in X_\tau$ (and thus for all
$z\in\C^g$):
\begin{equation}
\label{eq1}
\begin{array}{cl}
&\sum\limits_\e\T[\e](z)\T[\e](z)\T[\s](z)\hfill\\
=&\sum\limits_\e\sum\limits_{k=0}^{g}(-1)^{(\e+\s,e_{k+1})}\T[\e](z)\T[\e+s_k]
\T[\s+s_k](z),
\end{array}
\end{equation}
where we understand $e_{g+1}$ to be zero, and asumme moreover that
not all the coefficients appearing in front of $\T[\s+s_k](z)$ and
$\T[\s](z)$ are identically zero in $z$.
\end{thm}

To check that this makes sense let us do the computations in low
genus and see what we get. To simplify formulas, we write $[\e]$
for $\T[\e](z)$. We order the $\e$ for summation of the terms of
(\ref{eq1}) lexicographically to keep track of where we are.

\smallskip
\noindent {\bf Genus 2:} We do not expect to get any meaningful
equations, as any irreducible abelian variety of dimension two is
a hyperelliptic Jacobian, so our characterization should be
vacuous. We verify this; here is what formula (\ref{eq1}) yields
for $\s=00$:
$$
\begin{array}{cl}
\ &[00][00][00]+[01][01][00]+[10][10][00]+[11][11][00]\\
=&[00][00][00]+[01][01][00]-[10][10][00]-[11][11][00]\\
+&[00][10][10]-[01][11][10]+[10][00][10]-[11][01][10]\\
+&[00][11][11]+[01][10][11]+[10][01][11]+[11][00][11]
\end{array}
$$
and all the terms cancel. A similar computation shows that the
identity is also trivial for all other choices of $\s$.

\smallskip
\noindent {\bf Genus 3:} Here we have $\dim \J_3=\dim \A_3=6$,
while $\dim \I_3=5$, so we should have a non-trivial identity.
Indeed let us choose $\s=000$ and write down (\ref{eq1}) in this
case; after multiple cancellations and dividing by two it becomes
simply
$$
\begin{array}{cl}
&\T[000](z)\T[101](z)\T[101](z)+\T[011](z)\T[101](z)\T[110](z)\\
=&\T[010](z)\T[101](z)\T[111](z)+\T[001](z)\T[100](z)\T[101](z),
\end{array}
$$
which using formula (\ref{tT}) is equivalent to
\begin{equation}
\label{gen3}
\T[101](z)\cdot\tt{101}{111}(2z)\cdot\tt{101}{111}(0)=0.
\end{equation}

Choosing a different $\s$ yields a different equation: in general
we get $\T[101+\s](z)\cdot\tt{101}{111}(2z)\cdot
\tt{101}{111}(0)$. All of these equations together are equivalent
to $\tt{101}{111}(0)=0$, since theta functions of the second order
never all vanish simultaneously and $\tt{101}{111}(2z)$ cannot be
identically zero in $z$. Now to actually characterize
$\I_3\subset\A_3$ we need to get rid of the condition ``that the
basis of cycles is chosen as above'' in theorem \ref{thmchar}. But
choosing a different basis of cycles means acting on the period
matrix by a symplectic transformation. Since theta constants are
modular with respect to $\Gamma(2,4)$, conjugating equation
(\ref{gen3}) by any $\gamma\in\Gamma(2,4)$ would not change it, so
we only need to act by the finite group
$\operatorname{Sp}(2g,\Z)/\Gamma(2,4)$. It is well known that the
action of this group is transitive on the set of even theta
characteristics (see \cite{ig}), so we can get the vanishing of a
theta constant with any even characteristic.

\begin{prop} An irreducible abelian variety of genus 3 is a
hyperelliptic Jacobian if it has a vanishing theta constant with
even characteristic. This is known classically, see
\cite{mumford}.
\end{prop}

\smallskip
\noindent {\bf Genus 4:} here the situation is more interesting:
Mumford's conditions include some non-vanishing, so getting the
explicit equations for the closure of $Th(\I_4)\subset Th(\A_4)$,
without any inequalities serving to cut off the extra components
inside the reducible locus, would be interesting. For the case of
$\s=0000$ the equation we get from (\ref{eq1}) after cancellations
becomes the following cubic (we have rearranged the terms
lexicographically and omitted square brackets):
$$
\begin{array}{l}
0000\c1001\c1001+0000\c1010\c1010+0000\c1011\c1011+0000\c1101\c1101\\
0011\c1101\c1110+0101\c1000\c1101+0101\c1011\c1110+0110\c1010\c1100\\
0111\c1001\c1110+0111\c1011\c1100=0001\c1000\c1001+0001\c1100\c1101\\
0010\c1000\c1010+0010\c1101\c1111+0011\c1000\c1011+0100\c1010\c1110\\
0100\c1011\c1111+0101\c1001\c1100+0101\c1010\c1111+0110\c1001\c1111,
\end{array}
$$
while for example for $\s=0001$ we get
$$
\begin{array}{l}
0000\c1000\c1001+0000\c1100\c1101+0010\c1001\c1010+0011\c1001\c1011\\
0011\c1100\c1110+0100\c1000\c1101+0100\c1011\c1110+0101\c1010\c1110\\
0101\c1011\c1111+0111\c1000\c1110=0001\c1000\c1000+0001\c1010\c1010\\
0001\c1011\c1011+0001\c1100\c1100+0010\c1100\c1111+0100\c1001\c1100\\
0100\c1010\c1111+0110\c1000\c1111+0110\c1010\c1101+0111\c1011\c1101.
\end{array}
$$

Neither of these cubics is equal to $\T[\de](z)\,\tt\a\b(2z)
\tt\a\b(0)$ for any $\a,\b,\de$. However, from theorem
\ref{thmchar} we see that

\begin{prop} The vanishing of the full set of $16$ cubics similar
to the ones above, for all $\s$, identically in $z$ characterizes
a component of $Th(\I_4(2,4))\subset Th(\A_4^{irr}(2,4))$.
\end{prop}

By Mumford's and Poor's results such a component is also
determined by identical vanishing of some set of theta constants
with characteristics. Thus the vanishing of our 16 cubics should
imply, for irreducible abelian varieties, the vanishing of some
theta constants with characteristics and vice versa, but we are
now unable to see this directly.

The difficulty in doing so is not only due to the fact that the
cubic equations are very complicated. Indeed, thinking of each
cubic $f_\s$, evaluated at $z=0$, as a polynomial on $\P^{15}$
with zero locus $Z(f_\s)$, we can only say that $\mathop{\cap}
\limits_\s Z(f_\s)\cap Th(\A_4^{irr}(2,4))$ is contained in the
common zero locus of some quadrics (which are by (\ref{tT}) the
expressions for theta constants with characteristics in terms of
theta constants of the second order) on $Th(\A_4^{irr}(2,4))$. It
may in fact not be the case that the whole $\cap Z(f_\s)$ is
contained in the zero locus of these quadrics in $\P^{15}$. Thus
to be able to see the relation of the vanishing of our cubics to
the vanishing of theta constants with characteristics, we may need
to know the equations for the closure of
$Th(\A_4^{irr}(2,4))\subset \P^{15}$, which are not known.

The above discussion was for just one component of $Th(\I_4(2,4
))$, which projects to just one component of $\I_4\subset\A_4$.
The equations for the other components corresponding to different
choices of the basis of cycles on the curve are of course obtained
by acting on the set of 16 cubics by elements of
$\operatorname{Sp}(8,\Z)/\Gamma(2,4)$.

{\bf Final remark.} It seems likely that in any genus evaluating
equations (\ref{eq1}) only at $z=0$ for all $\s$ should yield the
defining set of equations for a component of $Th(\I_g(2,4))$.
Indeed using (\ref{change}) it can be easily shown that if these
are satisfied, then (\ref{eq1}) is satisfied for $z$ being any
point of order two. Thus both sides of (\ref{eq1}) are sections of
$6\T$ that agree at all points of order two, and one would hope
that then they agree everywhere and give the same function of $z$,
so that (\ref{eq1}) is true identically and we can apply theorem
\ref{thmchar}.

\end{document}